\documentclass[11pt,reqno,oneside]{amsart}
 \usepackage{latexsym}
 \usepackage{amssymb}
 \usepackage{amsfonts}
\usepackage{graphics}
\usepackage{mathcomp}
 \author[Mathonet]{P. Mathonet}
\thanks{University of Li\`ege, Institute of mathematics, Grande Traverse, 12 - B37, B-4000
 Li\`ege, Belgium email : P.Mathonet@ulg.ac.be, F.Radoux@ulg.ac.be}
\author[Radoux]{F. Radoux}
\date{\today} 
\title[Cartan connections and natural quantizations]{Cartan connections and natural and projectively equivariant quantizations}

\newtheorem{lem}{Lemma}
\newtheorem{thm}[lem]{Theorem}
\newtheorem{prop}[lem]{Proposition}

\theoremstyle{remark}
\newtheorem{rem}{Remark}
\theoremstyle{definition}
\newtheorem{defi}{Definition}

\newcommand{\R}{\mathbb{R}}
\newcommand{\Z}{\mathbb{Z}}
\newcommand{\C}{\mathbb{C}}

\newcommand{\N}{\mathbb{N}}

\renewcommand{\L}{\mathcal{L}}
\newcommand{\V}{\mathcal{V}}
\newcommand{\D}{\mathcal{D}}
\renewcommand{\S}{\mathcal{S}}

\newcommand{\g}{\mathfrak{g}}

\newcommand{\h}{\mathfrak{h}}

\newcommand{\euler}{\mathcal{E}}

\newcommand{\tree}[1][\gamma]{\mathcal{T}_{{#1}}}

\newcommand{\cc}{{\mathcal{C}}}
\begin{document}
\begin{abstract}
In this paper, we analyse the question of 
existence of a natural and projectively
equivariant symbol calculus, using the theory of projective Cartan
connections.
We establish a close relationship between the existence of such a natural 
symbol calculus and the existence of an
$sl(m+1,\R)$- equivariant calculus over $\R^m$ in the sense of
\cite{LO, BHMP}. Moreover we show that the formulae that hold in the
non-critical situation over $\R^m$ for
the $sl(m+1,\R)$- equivariant calculus can
be directly generalized to an arbitrary manifold by simply replacing the
partial derivatives by invariant differentiations with respect 
to a Cartan connection. 
\end{abstract}
%%%
\maketitle
%%%
\section{Introduction}
%%%
%%%
A quantization procedure can be
roughly defined as a linear bijection from the space of
classical observables to a space of differential operators acting on wave
functions (at least in the framework of geometric quantization, see
\cite{Woo}).

In our setting, the space $\mathcal{S}(M)$ of observables (also called
\emph{Symbols}) is made of smooth functions on
the cotangent bundle $T^*M$ of a manifold $M$, that are polynomial along the
fibres. The space $\mathcal{D}_{\lambda}(M)$ of differential operators consists
of operators acting on $\lambda$-densities over $M$.

It is known that
there is no natural quantization procedure. In other words, the spaces of
symbols and of differential operators are not isomorphic as representations of
$\mathrm{Diff}(M)$.

The idea of \emph{equivariant quantization}, introduced by P. Lecomte and
V. Ovsienko in \cite{LO} is to reduce the group of (local) diffeomorphisms
under consideration. 

They studied the case of the projective group
$PGL(m+1,\R)$ acting locally on the manifold $M=\R^m$ by linear fractional
transformations. 
They showed that the spaces of symbols and of differential operators are
canonically isomorphic as representations of $PGL(m+1,\R)$ (or its Lie algebra
$sl(m+1,\R)$). In other words, they showed that there exists a unique
\emph{projectively equivariant quantization}.

Independently of quantization purposes, the equivariant quantization map
proved to be a useful tool in the study of the spaces of differential
operators. Indeed, the inverse of such a map, called \emph{projectively
  equivariant symbol map}, is an equivariant bijection from the filtered space
$\mathcal{D}_{\lambda}(M)$ to its associated graded space $\mathcal{S}(M)$ and
allows to develop an equivariant symbol calculus.

In \cite{DO}, the authors studied the spaces $\mathcal{D}_{\lambda \mu}(\R^m)$ of
differential operators transforming $\lambda$-densities into $\mu$-densities
and their associated graded spaces $\mathcal{S}_{\delta}$.
They showed the existence and uniqueness of a projectively equivariant
quantization, provided the shift value $\delta=\mu-\lambda$ does not belong to
a set of critical values.

A first example of projectively equivariant symbol calculus for differential
operators acting on tensor fields was given in \cite{BHMP}.

At that point, all these results were dealing with a manifold endowed with a
 flat projective structure. In \cite{Bou1,Bou2}, S. Bouarroudj showed 
that the formula for the
 projectively equivariant quantization for differential operators of order two
 and three could be expressed using a torsion-free linear connection, in such
 a way that it only depends on the \emph{projective
 class} of the connection.

In \cite{Leconj}, P. Lecomte conjectured the existence of a quantization 
procedure 
\[Q:\S_{\delta}(M)\to{\mathcal{D}}_{\lambda\mu}(M)\]
depending on a torsion-free linear connection, that would be
 natural (in all arguments) and that would remain invariant under a projective
change of connection.

The existence of such a \emph{Natural and equivariant
   quantization procedure} was proved by M. Bordemann 
in \cite{Bor}, using the notion of Thomas-Whitehead connection  associated to
   a projective class of connections (see
   \cite{Thomas1,Whitehead,Roberts2,Roberts1,Roberts3} for a discussion of
   Thomas-Whitehead connections). 
His construction was later adapted by
   S. Hansoul in \cite{Hansoul} in order to deal with differential 
operators acting on forms, thus extending the results of \cite{BHMP}.

In \cite{MR}, we analysed the existence problem for a natural and projectively
 equivariant quantization using the theory of projective Cartan connections. 
We obtained an explicit formula for the natural and projectively equivariant
quantization in terms of the normal Cartan connection associated to a
projective class of linear connections. 

To our astonishment, it turned out that this explicit 
formula is nothing but the formula given in
\cite{DO} for the flat case, up to replacement of the partial derivatives by
invariant differentiations with respect to the Cartan connection.
 In particular, we
showed that the natural and projectively equivariant quantization exists if
and only if the $sl(m+1,\R)$-equivariant quantization exists in the flat case.

Recently, in her thesis, S. Hansoul \cite{sarah} showed how to generalize the method given
by M. Bordemann in order to solve the problem of existence of a natural and
projectively equivariant symbol calculus for differential operators acting on
arbitrary tensor fields.

In this paper, we analyse this question with the help of Cartan projective
connections. The advantage of our method is that it provides a direct
generalization of the formulae that can be written down in the context of
projectively equivariant symbol calculus over $\R^m$ (in the non-critical
situation), simply by replacing the
partial derivatives by invariant differentiations with respect to a Cartan
connection. It also provides a direct link between the existence of the
$sl(m+1,\R)$-equivariant quantization over $M=\R^m$ and the existence of a
natural and projectively equivariant quantization over an arbitrary manifold.
%%%%
\section{Problem setting}\label{tens}
In this section, we will describe the definitions of the spaces of
differential operators acting on tensor fields and of their corresponding
spaces of symbols. Then we will
set the problem of existence of projectively equivariant natural
symbol calculi.
\subsection{Tensor fields}
Let $M$ be a smooth manifold of dimension $m\geq 2$.
The arguments of the differential operators that we will 
consider are sections of
vector bundles associated to the linear frame bundle.
We consider irreducible representations of the group $GL(m,\R)$ defined as
follows :
let $(V,\rho_D)$ be the representation of $GL(m,\R)$ corresponding to a Young
diagram $Y_D$ of depth $n<m$. Fix $\lambda\in\R$ and $z\in \Z$ and set 
\[\rho(A)u=\vert det(A)\vert^{\lambda} (det(A))^z\rho_D(A)u,\]
for all $A\in GL(m,\R)$, and $ u\in V$.\\
If $(V,\rho)$ is such a representation, we denote by $V(M)$ the vector bundle
\[P^1M\times_{\rho}V.\]
We denote by $\V(M)$ the space of smooth sections of $V(M)$. 
 This space can be identified with the space $C^{\infty}(P^1M,V)_{GL(m,\R)}$ 
of functions $f$ such that
\[f(u A) = \rho (A^{-1}) f(u)\quad \forall u \in P^1M,\;\forall A\in
GL(m,\R).\]
Since $V(M)\to M$ is  associated to $P^1M$, there are natural
actions of $\mathrm{Diff}(M)$ and of $\mathrm{Vect}(M)$ on
$\V(M)$. 
%%%
\subsection{Differential operators and symbols}
%%%%
If $(V_1,\rho_1)$ and $(V_2,\rho_2)$ are representations of $GL(m,\R)$, we
 denote by $\mathcal{D}(\V_1(M),\V_2(M))$ (or simply by $\D(M)$ if there is no
 risk of confusion) the space of 
linear differential operators from $\V_1(M)$ to
$\V_2(M)$. The actions of $\mathrm{Vect}(M)$ and
$\mathrm{Diff}(M)$ are induced by their actions on $\V_1(M)$ and $\V_2(M)$ :
one has
\[(\phi\cdot D)(f) = \phi\cdot(D(\phi^{-1}\cdot f)),\quad\forall f\in \V_1(M)
, D\in \mathcal{D}(M),\mbox{and}\, \phi
\in \mathrm{Diff}(M).\]
The space $\mathcal{D}(M)$ is filtered by the order of
differential operators. We denote by $\mathcal{D}^k(M)$ the
space of differential operators of order at most $k$. It is well-known
that this filtration is preserved by the action of local diffeomorphisms.
The space of \emph{symbols}, which we will denote by $\S_{V_1,V_2}(M)$ or
simply by $\S(M)$, is then the graded space associated to $\mathcal{D}(M)$.

We denote by $S^l_{V_1,V_2}$ the vector space
 $S^l\R^m\otimes V_1^*\otimes V_2$. There is a natural representation
 $\rho$ of $GL(m,\R)$ on this space (the representation of $GL(m,\R)$ on
 symmetric tensors is the natural one). We then denote by 
$S^l_{V_1,V_2}(M)\to M$
 the
vector bundle
\[P^1M\times_{\rho}S^l_{V_1,V_2}\to M,\]
and by $\mathcal{S}^l_{V_1,V_2}(M)$ the space of smooth sections of
$S^l_{V_1,V_2}(M)\to M$, that is, the space 
$C^{\infty}(P^1M,S^l_{V_1,V_2})_{ GL(m,\R)}$.

 Then the \emph{principal symbol operator} $\sigma :
\mathcal{D}^l(M)\to \mathcal{S}^l_{V_1,V_2}(M)$ commutes with
the action of diffeomorphisms and is a bijection from the quotient space
$\mathcal{D}^l(M)/\mathcal{D}^{l-1}(M)$ to
$\mathcal{S}^l_{V_1,V_2}(M)$.
Hence the space of symbols is nothing but
\[\mathcal{S}(M)
=\bigoplus_{l=0}^{\infty}\mathcal{S}^l_{V_1,V_2}(M),\]
endowed with the classical actions of $\mathrm{Diff}(M)$ and of
$\mathrm{Vect}(M)$.
\subsection{Projective equivalence of connections}
We denote by $\mathcal{C}_M$ the space of torsion-free linear connections
on $M$. Two such connections are \emph{Projectively equivalent} if they define
the same geodesics up to parametrization. In algebraic terms, two connections
are projectively equivalent if there
exists a one form $\alpha$ on $M$ such that their associated covariant
derivatives $\nabla$ and $\nabla'$ fulfil the relation
\[\nabla'_XY = \nabla_XY + \alpha(X) Y + \alpha(Y)X.\]
This formulation was given by H. Weyl in \cite{Weyl}. 
%%%%%
\subsection{Natural and equivariant symbol calculus}
%%%%%
A (generalized) quantization on $M$ is a linear bijection $Q_{M}$
  from the space of symbols $\mathcal{S}_{V_1,V_2}(M)$ to the space 
of differential operators $\mathcal{D}(\V_1(M),\V_2(M))$ such
  that
\[\sigma(Q_{M}(T))=T,\;\forall T\in\mathcal{S}^{k}_{V_1,V_2}(M),\;\forall
  k\in\mathbb{N}.\]
The inverse of $Q_M$, denoted by $s_M$ is a symbol map.

Roughly speaking, a natural quantization map is a quantization map which 
depends on a
  torsion free linear connection $\nabla$  and commutes with the action of
  diffeomorphisms (on all arguments). More precisely, a natural quantization map is a collection of
  quantizations $Q_{M}$ (defined for every manifold $M$) depending on a
  torsion free linear connection on $M$ such that
\begin{itemize}
\item[(1)]
For all torsion free linear connections on $M$, $Q_{M}(\nabla)$ 
is a quantization map,
\item[(2)]
If $\phi$ is a local diffeomorphism from $M$ to $N$, then one has
\[Q_{M}(\phi^{*}\nabla)(\phi^{*}T)=\phi^{*}(Q_{N}(\nabla)(T)),\quad\forall 
\nabla\in\mathcal{C}(N),\,T\in\mathcal{S}(N)\]
 \end{itemize}
Such a natural quantization map is \emph{projectively equivariant} if one has
\[Q_{M}(\nabla)=Q_{M}(\nabla')\]
 whenever $\nabla$ and $\nabla'$ are projectively equivalent torsion free 
linear connections on $M$.  
%%%%
\section{The flat case}\label{flat}
%%%
Our goal in the present paper is to show the relationship between the
existence of a  natural and projectively
equivariant symbol calculus and the existence of an $sl(m+1,\R)$-equivariant
symbol calculus in the sense of \cite{LO}. 

In this section, we will briefly recall the problems addressed in the
 framework of projectively equivariant quantizations over vector spaces and
 outline the methods used in \cite{DLO} and \cite{BM}.
 However, we will generalize the tools of these papers and present them 
in an more algebraic and coordinate-free fashion. 
%%%
\subsection{Tensor fields, symbols and differential operators}
These objects were defined in section \ref{tens}, but when $M$ is the
Euclidean space $\R^m$, we make the following identifications :
\[\begin{array}{ccc}\mathcal{V}_{1}(\R^m)&\cong&C^{\infty}(\R^m,V_1),\\
\mathcal{S}^k_{V_1,V_2}(\R^m) & \cong & C^{\infty}(\R^m,S^k_{V_1,V_2}).
\end{array}\]
The Lie algebra $\mathit{Vect}(\R^m)$ acts on these spaces in a well-known
manner :  one has
\begin{equation}\label{der}(L_XT)(x)=X.T(x) - \rho_{*}(D_xX)T(x)\end{equation}
for every $X\in\mathit{Vect}(\R^m)$ and any symbol $T$.

The space $\mathcal{D}(\R^m)$ of differential operators 
is equipped with the Lie derivative $\L$ given by the commutator.
%%%%%%%%
\subsection{The projective algebra of vector fields}\label{sl}
%%%%%%%%%%
Consider the projective group $G=PGL(m+1,\R) =GL(m+1,\R)/\R_0\mbox{Id}$. Its Lie
algebra $\g=gl(m+1,\R)/\R\mbox{Id}$ is isomorphic to $sl(m+1,\R)$ and
decomposes into a direct sum of subalgebras
\[\g=\g_{-1}\oplus\g_0\oplus\g_1=\R^m\oplus gl(m,\R)\oplus\R^{m*}.\]
The isomorphism is given explicitly by 
\begin{equation}\label{psi}\psi : gl(m+1,\R)/\R\mbox{Id}\to\g_{-1}\oplus\g_0\oplus\g_1:\left[\left(\begin{array}{cc}A & v\\\xi & a
\end{array}\right)\right]\mapsto (v,A-a\,Id,\xi).\end{equation}
The group $G=PGL(m+1,\R)$ acts on $\R P^m$. Since $\R^m$ can be seen
 as the open
set of $\R P^m$ of equation $x^{m+1}=1$, there is a local action of $G$ on
$\R^m$. The vector fields associated to this action are given by 
\begin{equation}\label{real}\left\{\begin{array}{ccc}
X^h_x & = & -h\quad\mbox{if}\,h\in \g_{-1}\\
X^h_x & = & -[h,x]\quad\mbox{if}\,h\in \g_{0}\\
X^h_x & = & -\frac{1}{2}[[h,x],x]\quad\mbox{if}\,h\in\g_{1}\\
\end{array}\right.,\end{equation}
where $x\in\g_{-1}\cong\R^m$.
These vector fields define a subalgebra of $\mathrm{Vect}(\R^m)$, which is
isomorphic to $sl(m+1,\R)$.

Finally, it will be interesting for our computations to recall that the
subalgebra $\g_0$ is reductive and decomposes as
\begin{equation}\label{euler}\g_0=\h_0\oplus\R \euler\end{equation}
where $\h_0$ is (isomorphic to) $sl(m,\R)$ and where the \emph{Euler} element
$\euler$ is defined by $ad(\euler)\vert_{\g_{-1}}=-Id$. 
%%%%%%%%%%%%%
\subsection{Projectively equivariant quantizations}
%%%%%%%%%%%%
A projectively equivariant quantization (in the sense of \cite{DLO,LO}) is a
quantization
 \[Q :\mathcal{S}_{V_1,V_2}(\R^m)\to
\mathcal{D}(\R^m)\]
 such that, for every $X^h\in sl(m+1,\R)$, one has
\[\mathcal{L}_{X^h}\circ Q=Q\circ L_{X^h}.\]
The existence and uniqueness of such quantizations were discussed in
\cite{LO,DO,Lecras} for differential operators acting on densities and in
\cite{BHMP} for differential operators acting on forms.

Let us recall a first relationship, already mentioned in \cite{Lecras}, between natural and projectively equivariant
quantizations and $sl(m+1,\R)$-equivariant quantizations over $\R^m$ :
\begin{prop}
If $Q_{M}$ is a natural projectively equivariant quantization and if
$\nabla_0$ is the flat connection on $\R^m$, then
$Q_{\mathbb{R}^{m}}(\nabla_{0})$ is $sl(m+1,\R)$-equivariant.
\end{prop}

From now on to the end of this section, we will generalize the tools that were
used by Duval, Lecomte and Ovsienko in \cite{DLO}, and generalized in
\cite{BM} in order to 
obtain results of existence of such quantizations.
\subsection{The affine quantization map}
There exists a well-known bijection from symbols to differential operators
over $\R^m$ : the so-called \emph{Standard ordering} $Q_{\mathit{Aff}}$. 
If a symbol $T\in\mathcal{S}^k_{V_1,V_2}(\R^m)$ writes
\[T(x,\xi)=\sum_{\vert\alpha\vert=k}C_{\alpha}(x)\xi^{\alpha},\]
where $\alpha$ is a multi-index, $\xi\in\R^{m^*}$,  and $C_{\alpha}(x)\in V_1^*\otimes V_2$, 
then one has
\[Q_{\mathit{Aff}}(T)=\sum_{\vert\alpha\vert=k}C_{\alpha}(x)\circ(\frac{\partial}
{\partial x})^{\alpha}.\]
From our point of view, $Q_{\mathit{Aff}}$ is an affinely equivariant quantization
map. Indeed, it is easily seen that it exchanges the actions of the affine
algebra (made of constant and linear vector fields) on the space of symbols
and of differential operators. 

Now, we can use formula (\ref{real}) in order to express this quantization 
map in a coordinate-free manner :
\begin{prop}
If $h_1,\cdots, h_k\in\R^m\cong \g_{-1}$, $A\in V_1^*\otimes V_2$, $t\in
C^{\infty}(\R^m)$, and 
\[T(x)=t(x)\,A\otimes h_1\vee\cdots\vee h_k,\]
 one has
\[Q_{\mathit{Aff}}(T) = (-1)^k t(x)\circ A\circ L_{X^{h_1}}\circ\cdots\circ
L_{X^{h_k}}.\]
\end{prop}
\begin{proof}
The proof is straightforward. Simply notice that we are dealing with constant
vector fields.
\end{proof}
\subsection{The map $\gamma$}
Using the affine quantization map, one can endow the space of symbols with a
structure of representation of $\mathit{Vect}(\R^m)$, 
isomorphic to $\mathcal{D}(\R^m)$. Explicitly, we set
\[\mathcal{L}_XT=Q_{\mathit{Aff}}^{-1}\circ\mathcal{L}_X\circ 
Q_{\mathit{Aff}}(T),\]
for every $T\in\mathcal{S}(\R^m)$ and $X\in\mathit{Vect}(\R^m)$.

An equivariant quantization is then an $sl(m+1,\R)-$isomorphism from
the representation $(\S(\R^m), L)$ to the representation $(\S(\R^m),\L)$.

In order to measure the difference between these representations the map 
\[\gamma : \g\to gl(\S(\R^m),\S(\R^m)) : h\mapsto \gamma(h)=\L_{X^h}-L_{X^h}\]
was introduced in \cite{BM}. This map can be easily computed in coordinates,
and we recall here its most important properties :
\begin{prop}\label{gamma0}
The map $\gamma$ is a Chevalley-Eilenberg 1-cocycle and vanishes
  on $\g_{-1}\oplus \g_0$. 

For every $h\in\g_1$ and $k\in \N$, the restriction
  of $\gamma(h)$ to $\S^k(\R^m)$ has values in $\S^{k-1}(\R^m)$ 
and is a differential
  operator of order zero with constant coefficients.

For every $h,h'\in\g_1$, one has $[\gamma(h),\gamma(h')]=0$.
\end{prop}
{\bf Remark :} Since $\gamma(h)$ is a differential operator of order zero with
constant coefficients, it is completely determined by its restriction to
constant symbols.

For our purpose, it will be also interesting to obtain a coordinate free
expression of $\gamma$. We have 
\begin{prop}\label{gamma1}
For every $h_1,\cdots, h_k\in\R^m\cong \g_{-1}$, $A\in V_1^*\otimes V_2$
and $h\in \g_1\cong \R^{m*}$ we have
\[\begin{array}{r}\gamma(h)(h_1\vee\cdots\vee h_k\otimes A)=\sum_{i=1}^k
h_1\vee\cdots(i)\cdots\vee h_k\otimes (A\circ \rho_{1_*}([h_i,h]))\\
+\sum_{i=1}^k\sum_{j<i}h_1\vee\cdots(i,j)\cdots\vee
h_k\vee[h_j,[h_i,h]]\otimes A,
\end{array}\]
where $[h_i,h]$ is in $gl(m,\R)$ by the isomorphism (\ref{psi}).
\end{prop}
\begin{proof}
By the very definition of $\gamma$, the expression 
\[Q_{\mathit{Aff}}(\gamma(h)(h_1\vee\cdots\vee h_k\otimes A))\]
is equal to
\begin{equation}\label{gamma2}\mathcal{L}_{X^{h}}\circ Q_{\mathit{Aff}}(h_{1}\vee\cdots\vee h_{k}\otimes
A)- Q_{\mathit{Aff}}(L_{X^{h}}(h_{1}\vee\ldots\vee h_{k}\otimes A)).\end{equation}
This expression is a differential operator of order at most $k$. Its  term of
order $k$ vanishes : just apply the operator
$\sigma$ to (\ref{gamma2}). Hence, we only have to sum up the terms of
order less or equal to $k-1$ in the first term of expression (\ref{gamma2}).
This latter term writes
\[(-1)^{k}[L_{X^{h}}\circ A\circ L_{X^{h_{1}}}\circ\cdots\circ L_{X_{h_{k}}} -
  A\circ L_{X^{h_{1}}}\circ\cdots\circ L_{X^{h_{k}}}\circ L_{X^{h}}].\]
The first term is of order $k$ and $k+1$.
The second term is
\[\begin{array}{l} 
 -(-1)^k A\circ L_{X^{h}}\circ L_{X^{h_{1}}}\circ\cdots\circ L_{X^{h_{k}}}\\
-(-1)^{k}\sum_{i=1}^kA\circ L_{X^{[h_i,h]}}\circ
L_{X^{h_1}}\circ\cdots(i)\cdots \circ L_{X^{h_k}}\\
-(-1)^{k}\sum_{i=1}^k\sum_{j<i}A\circ  L_{X^{[h_j[h_i,h]]}}\circ L_{X^{h_1}}\circ\cdots(i,j)\cdots\circ L_{X^{h_k}},
\end{array}\]
since $[h_j[h_i,h]]$ is a constant vector field. \\The first summand is again of
order $k$ and $k+1$.
Now, in view of formula (\ref{der}), the term of order $k-1$ in 
\[-(-1)^{k}\sum_{i=1}^kA\circ L_{X^{[h_i,h]}}\circ
L_{X^{h_1}}\circ\cdots(i)\cdots \circ L_{X^{h_k}}\]
is exactly 
\[(-1)^{k}\sum_{i=1}^kA\circ \rho_{1_*}(DX^{[h_i,h]})\circ
L_{X^{h_1}}\circ\cdots(i)\cdots \circ L_{X^{h_k}},\]
and the result follows since $DX^{[h_i,h]}=-ad([h_i,h]).$
\end{proof}
Now, the map $\gamma$ has an important invariance property :
\begin{prop}\label{gogamma}
For every $a\in GL(m,\R)$, $h\in \g_1$ and every symbol $T$, one has
\begin{equation}\label{invgamma}
\rho(a)(\gamma(h)T)=\gamma(Ad(a)h)(\rho(a)T)
\end{equation}
\end{prop}
\begin{proof}
The proof is a straightforward computation, using proposition \ref{gamma1}.
\end{proof}
%%%%
\subsection{Casimir operators}
In \cite{DLO,BM}, the construction of the quantization is based on the
comparison of the spectra and of the eigenvectors of some (second order)
Casimir operators. These operators are on the one hand the Casimir operator
$C$ associated to the representation $(\S(\R^m), L)$ and on the other hand 
the Casimir operator $\cc$ associated to the representation
$(\S(\R^m),\L)$.
We will briefly adapt the results of \cite{BM} in order to compute
these operators.

From now on to the end of this section, we choose a basis
$(e_r,h_s,\euler,\epsilon^t)$ of $sl(m+1,\R)$, such that the bases $e_i$ and
$\epsilon^j$ of $\g_{-1}$ and $\g_1$ are Killing-dual and $h_j$ is a basis of
$\h_0$. It was then proved in \cite{BM} that the dual basis writes $(\epsilon^r,h_s^*,\frac{1}{2m}\euler,e_t)$
and that moreover there holds
\begin{equation}\label{crochet}
\sum_{r=1}^m[e_r,\epsilon^r]=-\frac{1}{2}\euler.
\end{equation}
We also set 
\[N=2\sum_i\gamma(\epsilon^i)L_{X^{e_i}}.\]
The following result is then the direct generalization of the corresponding
one in \cite{BM} :
\begin{prop}
The Casimir operators are related by
\begin{equation}\label{casinil}\cc=C+N.\end{equation}
\end{prop}
The next step is to analyse the eigenvalue problem for the operator $C$. To
this aim, we have to fix some more notation :  as a 
representation of $\h_{0}\cong sl(m,\R)$,
$S^k_{V_1,V_2}=S^k\R^m\otimes V_1^*\otimes V_2$ decomposes as a sum of
irreducible representations
\begin{equation}\label{dec}S^k_{V_1,V_2}=\oplus_{s=1}^{n_{k}}I_{k,s}.
\end{equation}
For each irreducible representation $I_{k,s}$ we denote by 
$E_{k,s}$ the corresponding space of sections, that is 
\[E_{k,s} = C^{\infty}(\R^m,I_{k,s}).\]
Furthermore, in $sl(m,\C)$, we consider the usual Cartan subalgebra
$\mathfrak{C}$ made of diagonal and traceless matrices. We consider the elements
of $\mathfrak{C}^*$ defined by 
\[\delta_i(diag(a_1,\cdots,a_m))=a_i.\]
It is known that a simple root system is given by $\{\delta_i-\delta_{i+1},
(i=1,\ldots,m-1)\}$.
The Weyl vector is defined as half the sum of the positive roots and is given
by 
\[\rho_{S}=\sum_i (m-i)\delta_i.\]
The Killing form of $sl(m,\C)$ is
the extension of the Killing form of $sl(m,\R)$ and induces a scalar product $(\,,\,)$
on the real vector space spanned by the roots. 

Now for each irreducible representation $I_{k,s}$ of $sl(m,\R)$, the
complexified representation $I_{k,s}\otimes \C$ of $sl(m,\C)$ is also
irreducible and we denote by $\mu_{I_{k,s}}$ its highest weight.\\
In view of the definitions of the representations $V_1$ and $V_2$, there exist
real numbers $a_1$ and $a_2$ such that
\[\rho_{*}(Id)\vert_{V_i}=a_i\,Id.\]
In order to be consistent we respect to the definition of the shift value in
\cite{BM,DLO}, we define the shift of the pair $(V_1,V_2)$ by
\[\delta=\frac{1}{m}(a_1-a_2).\] 
We can now state the main result concerning the operator $C$ :
\begin{thm}\label{casival}
The space of symbols $\S(\R^m)$ is the direct sum of eigenspaces of $C$.
Precisely, for every $k\in\N$, the restriction of 
  $C$ to $E_{k,s}$ is equal to
$\alpha_{k,s}\mbox{Id}_{E_{k,s}}$ where
\begin{equation}\label{eq.casival}
\alpha_{k,s}=\frac{1}{2m}(m\delta -k)(m(\delta -1)-k)
+\frac{m}{m+1}(\mu_{I_{k,s}}, \mu_{I_{k,s}}+ 2 \rho_{S}).
\end{equation}
\end{thm}
\begin{proof}
With our choice of dual bases, the operator $C$ writes
\[\sum_i( L_{X^{\epsilon^i}}\circ L_{X^{e_i}} + L_{X^{e_i}}\circ
L_{X^{\epsilon^i}})+\frac{1}{2m}(L_{X^{\euler}})^2 + \sum_j L_{X^{h_j}}\circ
L_{X^{h_j^*}},\]
that is,  using relation (\ref{crochet}), 
\[ 2 \sum_i( L_{X^{\epsilon^i}}\circ L_{X^{e_i}}) -\frac{1}{2}L_{X^{\euler}}+\frac{1}{2m}(L_{X^{\euler}})^2 + \sum_j L_{X^{h_j}}\circ
L_{X^{h_j^*}}.\]
Since $C$ commutes with $L_{X^h}$ for all $h\in \g_{-1}\cong\R^m$ it must have constant
coefficients. Hence, we only take such terms into account. We then use the
expression of the Lie derivative (\ref{der}), and the expression of the
vector fields $X^h$, and get
\[C=-\frac{1}{2}\rho_*(ad(\euler))+\frac{1}{2m}(\rho_*(ad(\euler)))^2 + \sum_j
\rho_*(ad(h_j))\circ \rho_*(ad(h_j)),\]
where $ad$ is the adjoint representation of $\g_0$ on $\g_{-1}$.
The terms involving the Euler element can be easily computed since
$ad(\euler)\vert_{\g_{-1}}=-Id$. Now, the restriction of the Killing form of $sl(m+1,\R)$ to
the subalgebra $sl(m,\R)$ is $\frac{m+1}{m}$ times the Killing form of
$sl(m,\R)$, so that the bases $(h_j)$ and $(\frac{m+1}{m}h_j^*)$ are dual with
respect to the latter. Then
\[ \sum_j
\rho_*(ad(h_j))\circ \rho_*(ad(h_j))=\frac{m}{m+1}C',\]
where $C'$ is the Casimir operator of $sl(m,\R)$ acting on $I_{k,s}$ or the
Casimir operator of $sl(m,\C)$ acting on $I_{k,s}\otimes \C$. It is well known 
(see for instance~\cite[p.~122]{hum}) that this operator equals
\[(\mu_{I_{k,s}}, \mu_{I_{k,s}}+ 2 \rho_{S})\]
times the identity of $I_{k,s}$.
\end{proof}
%%%
\subsection{Tree-like subspaces and critical situations}
%%%
In order to analyse the spectrum of the operator $\cc$ we introduce, as in
\cite{BM,DLO}, the treelike susbspaces associated to an irreducible
representation $I_{k,s}\subset S^k_{V_1,V_2}$ : we set
\[
\tree(I_{k,s})=\bigoplus_{l\in\N}\tree^l(I_{k,s}),
\]
where $\tree^0(I_{k,s})=I_{k,s}$ and
$\tree^{l+1}(I_{k,s})=\gamma(\g_1)(\tree^l(I_{k,s}))$,
for all $l\in\N$.
The spaces $\tree^l(E_{k,s})$ are defined in the same way, so that 
\[\tree^l(E_{k,s})= C^{\infty}(\R^m,\tree^l(I_{k,s}))\]
for all $l\in\N$.

The spaces $\tree(E_{k,s})$ have a nice property :
\begin{prop}\label{stable}
The space $\tree(E_{k,s})$ is stable under the actions of $L_{X^h}$ and
$\L_{X^h}$, for every $h\in sl(m+1,\R)$.
\end{prop}
\begin{proof}
Proposition \ref{gogamma} allows to prove by induction that $\tree(E_{k,s})$
is stable under the action of $\rho_*(A)$ for all $A\in gl(m,\R)$. It is then
clearly stable under the Lie derivative $L_{X^h}$ for all $h\in sl(m+1,\R)$
because of the expression (\ref{der}) of $L_{X^h}$. The result then follows
since $\L_{X^h}=L_{X^h}+\gamma(h).$
\end{proof}
The following definition is a direct generalization of the ones of \cite{BM,
  DLO} :
\begin{defi}
An ordered pair of representations $(V_1,V_2)$
is \emph{critical} 
if there exists $k,s$ such that the eigenvalue $\alpha_{k,s}$ 
belongs to the spectrum of the restriction of $C$ to 
$\bigoplus_{l\geq 1}\tree^l(E_{k,s})$.
\end{defi}
%%%%
\subsection{Construction of the quantization}\label{flatconst}
%%%%%
The result is the following :
\begin{thm}\label{flatex}
If $(V_1,V_2)$ is not critical, then there exists a
projectively equivariant quantization from $\mathcal{S}_{V_1,V_2}(\R^m)$ to $\mathcal{D}(\V_1(\R^m),\V_2(\R^m))$.
\end{thm}
\begin{proof}
The proof is as in \cite{BM} and \cite{DLO}. We give here the main ideas for
the sake of completeness.

First remark that for every $T\in E_{k,s}$ there exists a unique eigenvector
$\hat{T}$ of $\cc$ with eigenvalue $\alpha_{k,s}$ such that 
\[\left\{\begin{array}{l}\hat{T}= T_k + T_{k-1}+\cdots+T_0,\quad T_k=T\\
T_l\in \tree^{k-l}(E_{k,s})\quad\mbox{for all}\, l\leq
k-1.\end{array}\right.\]
Indeed, these conditions write 
\begin{equation}\label{flatP}\left\{\begin{array}{l}
C(T)=\alpha_{k,s}T\\
(C-\alpha_{k,s}\mbox{Id}))T_{k-l}=- N(T_{k-l+1})\quad\forall
l\in\{1,\cdots,k\}\\
T_{k-l}\in\tree^{l}(E_{k,s})
\end{array}\right.\end{equation}
The first condition is satisfied since $T$  belongs to $E_{k,s}$. For the
second and third ones, remark that if $T_{k-l+1}$ is in $\tree^{l-1}(E_{k,s})$
then  $N(T_{k-l+1})$ belongs to  $\tree^{l}(E_{k,s})$, by proposition
\ref{stable}. 
Now, $\tree^{l}(E_{k,s})$ decomposes as a direct sum of eigenspaces of $C$, as
theorem \ref{casival} shows. The restriction of the operator
$C-\alpha_{k,s}\mbox{Id}$ to each of these eigenspaces is a non-vanishing 
scalar multiple of the identity,
hence the existence and uniqueness of $T_{k-l}$.

Now, define the quantization $Q$ by
\[Q\vert_{E_{k,s}}(T)=\hat{T}.\]
It is clearly a bijection. 

It also fulfils 
\[Q\circ L_{X^h} = \L_{X^h}\circ Q\quad\forall h\in\,sl(m+1,\R).\]
Indeed, for all $T\in E_{k,s}$, $Q(L_{X^h}T)$ and $ \L_{X^h}(Q(T))$ share the
following properties
\begin{itemize}
\item They are eigenvectors of $\cc$ of eigenvalue $\alpha_{k,s}$ because on
  the one hand $\cc$ commutes with $\L_{X^h}$ for all $h$ and on the other
  hand $L_{X^h}T$ belongs to $E_{k,s}$ by proposition \ref{stable},
\item their term of degree $k$ is exactly $L_{X^h}T$,
\item they belong to $\tree(E_{k,s})$ by proposition \ref{stable}.
\end{itemize}
The first part of the proof shows that they have to coincide.
\end{proof}
\subsection{A technical result}
%%%
The following proposition will be fundamental for our purpose :
\begin{prop}\label{tech}
The relation 
\[ [\gamma(h),C]=2\sum_{i}\gamma(\epsilon^{i})\rho_{*}([h,e_{i}]),\]
holds for all $h\in\g_{1}.$
\end{prop}
\begin{proof}
As a Casimir operator commutes with the corresponding
  representation, we have :
\[\begin{array}{lll}
[\mathcal{L}_{X^{h}}, \mathcal{C}] &=&0,\\
  {[L_{X^{h}}, C]} & = & 0.
\end{array}\]
These equations lead to
\[[L_{X^{h}}, N]+[\gamma(h), C]+[\gamma(h), N]=0.\]
It is easily seen using proposition \ref{gamma0} that
$[\gamma(h), N]$ vanishes. Moreover, we have
\[[L_{X^{h}}, N] = 2\sum_{i}(
L_{X^{h}}\gamma(\epsilon^{i})L_{X^{e_i}}-
\gamma(\epsilon^{i})L_{X^{e_i}}L_{X^{h}}).\]
The terms of order greater than zero in this expression must vanish because
$[\gamma(h), C]$ has order zero. Hence we only collect the terms of order
zero.

Using formula (\ref{der}), we see that the first term is of order greater or
equal to 1.
The second term writes
\[-2\gamma(\epsilon^{i})\sum_{i}(L_{X^{h}}L_{X^{e_i}}+L_{X^{[e_i,h]}})\]
The terms of order zero in this expression are
\[2\sum_{i}\gamma(\epsilon^{i})\rho_{*}(DX^{[e_i,h]}),\]
hence the result.
\end{proof}
%%%%%
%%%%%
\section{Tools of the curved case}
%%%%%
Here we will adapt the tools presented in section \ref{flat} 
to the curved situation.
\subsection{Projective structures and Cartan projective connections}
These tools were presented in detail in \cite[Section 3]{MR}. We give here the
 most important ones for this paper to be self-contained.

 We consider the group $G=PGL(m+1,\R)$ described in section \ref{sl}.
 We denote by $H$ 
the subgroup associated to the subalgebra $\g_0\oplus\g_1$, that is 
\begin{equation}\label{hproj}H=\{\left(\begin{array}{cc}A & 0\\\xi & a
\end{array}\right) : A\in GL(m,\R),\xi\in \R^{m*}, a\not = 0\}/\R_0\mbox{Id},
\end{equation}
The group $H$ is the semi-direct product $G_0 \rtimes G_1$, where
$G_0$ is isomorphic to $GL(m,\R)$ and $G_1$ is isomorphic to $\R^{m*}$. 

Then there is a projection 
 \[\pi : H\to GL(m,\R) : \left[\left(\begin{array}{cc}A & 0\\\xi & a
\end{array}\right)\right]\mapsto \frac{A}{a}\]
It is well-known that $H$ can be seen as a subgroup of the group of 2-jets 
$G_m^2$, an explicit formula is given in \cite{Kobabook, MR}.\\
A \emph{Projective
  structure on $M$} is then a reduction of the second order jet-bundle $P^2M$
to the group $H$. \\
The following result (\cite[P. 147]{Kobabook}) is 
the starting point of our method :
\begin{prop}[Kobayashi-Nagano]
There is a natural one to one correspondence between the projective
 equivalence classes of torsion-free linear connections on $M$
and the projective structures on $M$.
\end{prop}
We now recall the definition of a projective Cartan connection :
\begin{defi}
Let $P\to M$ be a principal $H$-bundle. A projective Cartan connection on $P$
is an $sl(m+1,\R)$- valued 1-form $\omega$ such that
\begin{itemize}
\item There holds $R_a^*\omega=Ad(a^{-1})\omega,\quad\forall a \in H$,
\item One has $\omega(k^*)=k\quad\forall k\in h=\g_0\oplus \g_1$,
\item For all $u\in P$, $\omega_u :T_uP\to sl(m+1,\R)$ is a linear bijection.
\end{itemize}
\end{defi}
In general, if $\omega$ is a Cartan connection defined on a $H$-principal bundle $P$, then its
curvature $\Omega$ is defined as usual by
 \begin{equation}\label{curv} 
\Omega = d\omega+\frac{1}{2}[\omega,\omega].
\end{equation}
The notion of \emph{Normal} Cartan connection is defined by natural conditions
imposed on the components of the curvature.

Now, the following result (\cite[p. 135]{Kobabook}) gives the relationship between  projective
structures and Cartan connections :
\begin{prop}
 A unique normal
 Cartan projective connection is
 associated to every projective structure $P$. This association is natural.
\end{prop}
The connection associated to a projective structure $P$ is called the normal
projective connection of the projective structure.
%%%%
\subsection{Lift of equivariant functions}\label{Lift}
%%%%
As we continue, we will need to know the
relationship between equivariant functions on $P^1M$ and equivariant functions on
$P$. The following results were already quoted in \cite[p. 47]{Capinv}.

If $(V,\rho)$ is a representation of $\mathit{GL}(m,\R)$, then we define
a representation $(V,\rho')$ of $H$ by
\[\rho': H\to GL(V): \left[\left(\begin{array}{cc}A &0\\ \xi &
    a\end{array}\right)\right]\mapsto  \rho\circ \pi( \left[\left(\begin{array}{cc}A &0\\ \xi &
    a\end{array}\right)\right])
= \rho(\frac{A}{a})\]
for every $ A\in \mathit{GL}(m,\R), \xi\in
\R^{m*}, a\not=0$.

Now, using the representation $\rho'$,  we can give the relationship between
 equivariant functions on $P^1M$ and equivariant functions on $P$ :
If $P$ is  a projective structure on $M$, the natural projection $P^2M\to
P^1M$ induces a projection $p :P\to P^1M$ and we have a well-know result:
\begin{prop} If $(V,\rho)$ is a representation of $GL(m,\R)$, then the map
\[p^* : C^{\infty}(P^1M,V)\to C^{\infty}(P,V) : f \mapsto f\circ p\]
defines a bijection from $C^{\infty}(P^1M,V)_{\mathrm{GL}(m,\R)}$ to
$C^{\infty}(P,V)_{H}$.
\end{prop}
Now, since $\R^m$ and $\R^{m*}$ are natural representations of $GL(m,\R)$,
they become representations of $H$ and we can state an important property of the invariant differentiation :
\begin{prop}\label{gonabla}If $f$ belongs to $C^{\infty}(P,V)_{G_0}$ then
  $\nabla^{\omega}f\in C^{\infty}(P,\R^{m*}\otimes V)_{G_0}$.
\end{prop}
\begin{proof}The result is a direct consequence of the Ad-invariance of the
  Cartan connection $\omega$. 
\end{proof}
The main point that we will discuss in the next sections is that this result
is not true in general for $H$-equivariant functions : for an $H$-equivariant
function $f$, the function
$\nabla^{\omega}f$ is in general not $G_1$-equivariant.

As we continue, we will use the representation $\rho'_*$ of the Lie algebra of
$H$ on $V$. If we recall that this algebra is isomorphic to
 $gl(m,\R)\oplus \R^{m*}$
then we have 
\begin{equation}\label{rho}\rho'_* (A, \xi) = \rho_*(A),\quad\forall A\in
  gl(m,\R), \xi\in \R^{m*}.\end{equation}
In our computations, we will make use of the infinitesimal version of the
equivariance relation : If $f\in C^{\infty}(P,V)_H$ then one has
\begin{equation}\label{Invalg}
L_{h^*}f(u) + \rho'_*(h)f(u)=0,\quad\forall h\in gl(m,\R)\oplus\R^{m*}\subset sl(m+1,\R),
\forall u\in P.
\end{equation}
\subsection{The curved affine quantization map}
The construction of the curved analog of the affine quantization map is based
on the concept of invariant differentiation developed in
\cite{Capinv,Capcart}.
Let us recall the definition :
\begin{defi} Let $(V, \rho)$ be a representation of $H$. If $f\in
 C^{\infty}(P,V)$, then the invariant differential of $f$ with respect to
$\omega$ is the function \\
$\nabla^{\omega}f\in C^{\infty}(P,\R^{m*}\otimes V)$
defined by
\[ \nabla^{\omega}f(u)(X) = L_{\omega^{-1}(X)}f(u)\quad\forall u\in
P,\quad\forall X\in\R^m.\]
\end{defi}
We will also use an iterated and symmetrized version of the invariant 
differentiation 
\begin{defi}
If $f\in C^{\infty}(P,V)$ then $(\nabla^{\omega})^k f
\in C^{\infty}(P,S^k\R^{m*}\otimes V)$ is defined by
\[(\nabla^{\omega})^k f(u)(X_1,\ldots,X_k) = \frac{1}{k!}\sum_{\nu}
L_{\omega^{-1}(X_{\nu_1})}\circ\ldots\circ L_{\omega^{-1}(X_{\nu_k})}f(u)\]
for $X_1,\ldots,X_k\in\R^m$.
\end{defi}
Using this iterated version, we can transform a symbol $T\in C^{\infty}
(P,S^k_{V_1,V_2})$ into a differential operator $Q_{\omega}(T)$ acting on
functions $f\in C^{\infty}(P,V_1)$ by setting
\begin{equation}\label{Qom}
Q_{\omega}(T)(f)=\langle T,(\nabla^{\omega})^k f\rangle.
\end{equation}
%%%%%
Explicitly, when the symbol $T$ writes $t A\otimes h_1\vee\cdots\vee h_k$ 
for $t\in C^{\infty}(P)$, $A\in V_1^*\otimes V_2$ and $h_1,\cdots,
h_k\in\R^m\cong\g_{-1}$ then one has
\[Q_{\omega}(T)f=\frac{1}{k!}\sum_{\nu}t A\circ
L_{\omega^{-1}(h_{\nu_1})}\circ\cdots\circ L_{\omega^{-1}(h_{\nu_k})}f,\]
where $\nu$ runs over all permutations of the indices $\{1,\cdots,k\}$ and $t$
is considered as a multiplication operator\vspace{0,2cm}.

\begin{rem}\label{remf} If $T\in C^{\infty} 
(P,S^k_{V_1,V_2})$ is
$H-$equivariant, the differential operator $Q_{\omega}(T)$ does not transform
$H-$equivariant functions into $H-$equivariant functions. Indeed, when $f$ is
$H-$equivariant, the function $(\nabla^{\omega})^k f$ is only
$G_0-$equivariant. Hence the function $Q_{\omega}(T)f$ does not correspond to
a section of $\mathcal{V}_2(M)$. As we continue, we will show that one can modify the
symbol $T$ by lower degree correcting terms in order to solve this problem.
\end{rem}
\subsection{Measuring the default of equivariance}
Throughout this section, $T$ will denote an element of $C^{\infty}
(P,S^k_{V_1,V_2})_{G_0}$ and $f\in C^{\infty}(P,V_1)_{G_0}$ 
(remark that this ensures that
$Q_{\omega}(T)(f)$ is in $C^{\infty}(P,V_2)_{G_0}$). 
Now, in order to analyse the invariance of functions, we have this first easy
result
\begin{prop}If $(V,\rho)$ is a representation of $G_0$ and becomes a
representation of $H$ as stated in section \ref{Lift}, then a function $v\in
 C^{\infty}(P,V)$ is $H-$equivariant iff 
\[\left\{\begin{array}{l}
\mbox{$v$ is $G_0-$equivariant}\\
\mbox{One has $L_{h^*}v=0$ for every $h$ in $\g_1$}
\end{array}\right.\]
\end{prop}
\begin{proof}
Of course, $H-$equivariance is equivalent to $G_0-$ and
$G_1-$equivariance. Now, $G_1-$equivariance is equivalent to
$\g_1-$equivariance since $G_1$ is a vector space. The result follows since
$G_1$ acts trivially on $V$.
\end{proof}
Since basically, our tools preserve the $G_0$-equivariance, we are mostly
interested in the $\g_1$-equivariance.
The following result is the keystone of our method :
\begin{prop}
The relation
\[L_{h^{*}}Q_{\omega}(T)(f)-Q_{\omega}(T)(L_{h^{*}}f)=Q_{\omega}
((L_{h^*}+\gamma(h))T)(f)\]
holds for all $f\in C^{\infty}(P,V_1)_{G_0}$,
$h\in\mathfrak{g}_{1}$, and $T\in C^{\infty}(P,S^{k}_{V_1,V_2})$.
\end{prop}
\begin{proof}
The relation that we have to prove writes
$$\langle T, L_{h^*}(\nabla^{\omega})^kf-(\nabla^{\omega
  })^kL_{h^*}f\rangle=\langle\gamma(h)T,(\nabla^{\omega})^{k-1}  f\rangle.$$
Since both sides are $C^{\infty}(P)$-linear in $T$, it is sufficient to check
  this relation for a constant symbol $T$ that has the form
  $X^k\otimes A$, where $X\in\mathfrak{g}_{-1}$ and
  $A\in V_1^*\otimes V_2$. Then the left-hand side writes
$$A(L_{h*}L_{\omega^{-1}(X)}\ldots L_{\omega^{-1}(X)}f-L_{\omega^{-1}(X)}\ldots
L_{\omega^{-1}(X)}L_{h^*}f)$$
and is equal to
\[\begin{array}{l}A(\sum_{j=1}^{k} L_{\omega^{-1}(X)}\ldots\overset{(j)}{L_{[h,X]^{*}}}\ldots L_{\omega^{-1}(X)}f)\\
=A(\sum_{j=1}^{k-1}\sum_{i>j}L_{\omega^{-1}(X)}\overset{(j)}{\ldots}
\overset{(i)}{L_{\omega^{-1}([[h,X],X])}}\ldots L_{\omega^{-1}(X)}f)\\
-A(\sum_{i=1}^{k}L_{\omega^{-1}(X)}\overset{(i)}{\ldots}
L_{\omega^{-1}(X)}\rho_{1_*}([h,X])f).\end{array}\]
Moreover, one has 
\[ [[h,X],X]=-2\langle h,X\rangle X\]
and the result follows by proposition \ref{gamma1}.
\end{proof}
%%%%%
\subsection{Curved Casimir operators}
%%%%
The parallelism between the flat and curved situation suggests to define an
analog of the operator $\mathcal{C}$. 

We first define the analog of $N$ by setting
$$N^{\omega}=-2\sum_{i}\gamma(\epsilon^{i})L_{\omega^{-1}(e_i)}.$$
Then we can define the operators $C^{\omega}$ and $\mathcal{C}^{\omega}$ by
their restrictions to the spaces $C^{\infty}(P,I_{k,s})$ :
for all $T\in C^{\infty}(P,I_{k,s})$, we set
%%%
\[\left\{\begin{array}{lll}C^{\omega}(T)&=&\alpha_{k,s}T\\

\mathcal{C}^{\omega}(T)&=&C^{\omega}(T)+N^{\omega}(T),
\end{array}\right.\]
where $\alpha_{k,s}$ is the eigenvalue of $C$ on $E_{k,s}=C^{\infty}(\R^m,I_{k,s})$.

The operator $\cc^{\omega}$ has the following property
\begin{prop}\label{commute}
For every $h\in\mathfrak{g}_{1}$, one has 
$$[L_{h^{*}}+\gamma(h), \mathcal{C}^{\omega}]=0$$
on $C^{\infty}(P,S^k_{V_1,V_2})_{G_0}$.
\end{prop}
\begin{proof}
First we have
\[ [L_{h^{*}}+\gamma(h), C^{\omega}+N^{\omega}]= [L_{h^{*}},
  C^{\omega}]+[L_{h^{*}}, N^{\omega}]+[\gamma(h), C^{\omega}]+[\gamma(h),
  N^{\omega}].\]
Next $[L_{h^{*}},C^{\omega}]=0$ because $L_{h^{*}}$ stabilizes each
  eigenspace $C^{\infty}(P,I_{k,s})$ of $C^{\omega}$. In the same way, we
  have
$[\gamma(h),N^{\omega}]=0$ since 
\begin{itemize}
\item by proposition \ref{gamma0}, $[\gamma(h),\gamma(\epsilon^i)]=0,$ 
\item $[\gamma(h), L_{\omega^{-1}(e_i)}]=0$ because $\gamma(h)$ only acts 
on the
  target space $S^k_{V_1,V_2}$.
\end{itemize}
Eventually, we have
\[\begin{array}{lll} [L_{h^{*}}, N^{\omega}]&=& - 2 \sum_i
\gamma(\epsilon^i)[L_{h^*},L_{\omega^{-1}(e_i)}]\\
&=&-2\sum_i\gamma(\epsilon^i)L_{[h,e_i]^*}\\
&=&2\sum_{i}\gamma(\epsilon^{i})\rho_{*}([h,e_{i}]),\end{array}\]
by the $G_{0}$-equivariance. We conclude by proposition \ref{tech}.
\end{proof}
For the operator $N^{\omega}$, we have the following result :
\begin{prop}\label{Goinv}
The operator $N^{\omega}$ preserves the $G_0$-equivariance of functions.
\end{prop}
\begin{proof}
This property is a consequence of the proposition \ref{gogamma} and of the fact that the invariant differentiation
preserves the $G_{0}$-equivariance. One has successively, for all $f\in C^{\infty}(P,S^k_{V_1,V_2})$, $u\in P$ and $g\in G_{0}$ :
\[\begin{array}{lll}
(N^{\omega}(f))(ug)&=&
  \sum_{i}(\gamma(\epsilon^{i})L_{\omega^{-1}(e_i)}f)(ug)\\ & = &
  \sum_{i}\gamma(\epsilon^{i})(\nabla^{\omega}f)(ug)(e_i) \\
                                                        & = &
  \sum_{i}\gamma(\epsilon^{i})\rho(g^{-1})((\nabla^{\omega}f)(u)(Ad(g)e_i)) \\
                                                        & = &
  \sum_{i}\rho(g^{-1})(\gamma(Ad(g)\epsilon^{i})
(\nabla^{\omega}f)(u)(Ad(g)e_i))
  \\
                                                        & = &
  \rho(g^{-1})\sum_{i}\gamma(\epsilon^{i})L_{\omega^{-1}(e_i)}f(u), \\
\end{array}\]
hence the result.
\end{proof}
%%%%
\section{Construction of the quantization}
The construction of the quantization is based on the eigenvalue problem for
the operator $\mathcal{C}^{\omega}$.

First remark that the construction of section \ref{flatconst} is still valid
in the curved case.
\begin{thm}\label{hat}
If the pair $(V_1,V_2)$ is not critical, for every $T$
in $C^{\infty}(P,I_{k,s})$, there exists a unique function $\hat{T}$ in $ C^{\infty}(P,\tree(I_{k,s}))$
 such that
\begin{equation}\label{curvP}\left\{\begin{array}{lll}
\hat{T}&=&T_k+\cdots+T_0,\quad T_k=T\\
\mathcal{C}^{\omega}(\hat{T})&=&\alpha_{k,s}\hat{T}.
\end{array}\right.
\end{equation} 
Moreover, if $T$ is $G_0$-invariant, then $\hat{T}$ is $G_0$-invariant.
\end{thm}
\begin{proof}
The function $\hat{T}$ exists and is unique.
Simply notice that the conditions in (\ref{curvP}) are nothing but the ones in
(\ref{flatP}) with $C$ replaced by $C^{\omega}$ and $N$ replaced by
$N^{\omega}$. Now, the main point in order to solve (\ref{flatP}) was that $N$
maps $\tree^{l-1}(E_{k,s})$ into $\tree^{l}(E_{k,s})$. This latter fact is
actually a property of $\gamma$ and therefore we have 
\[N^{\omega}(C^{\infty}(P,\tree^{l-1}(I_{k,s})))\subset C^{\infty}(P,\tree^{l}(I_{k,s})).\]
Suppose in addition that $T$ is $G_0$-invariant, then $\hat{T}$ is $G_0$-invariant. Indeed, $\hat{T}$ is obtained from
$T$ by applying successively the operators $N^{\omega}$ and the projectors
from $\tree^{l}(I_{k,s}))$ onto its irreducible components, and these operations preserve the $G_0$-equivariance of functions (see prop \ref{Goinv}).
\end{proof}
%%%%
This result allows to define the main ingredient in order to define the
quantization.
\begin{defi}
Suppose that the pair $(V_1,V_2)$ is not critical. Then
the map 
\[Q : C^{\infty}(P,S_{V_1,V_2})\to  C^{\infty}(P,S_{V_1,V_2})\]
is the linear extension of the association $T\mapsto \hat{T}$.
\end{defi}
The map $Q$ has the following nice property :
\begin{prop}
There holds
\begin{equation}
\label{Q}(L_{h^*}+\gamma(h))Q(T)=Q(L_{h^*}T),
\end{equation}
for every $h\in\g_1$ and every $T\in C^{\infty}(P,S_{V_1,V_2})_{G_0}$.
\end{prop}
\begin{proof}
The proof is just an adaptation of the one of theorem \ref{flatex}. 
It is sufficient to check the property for $T\in C^{\infty}(P,I_{k,s})_{G_0}$ 
(for all
$k$ and $s$). For such a $T$, the function $Q(L_{h^*}T)$ is defined by 
(\ref{curvP}) : it is the unique eigenvector
of $\cc^{\omega}$ in $C^{\infty}(P,\tree(I_{k,s}))$ corresponding to
the eigenvalue $\alpha_{k,s}$ and with highest degree term $L_{h^*}T$.
The left hand side of equation (\ref{Q}) has $L_{h^*}T$ as highest degree term
since $\gamma(h)$ lowers the degree of its arguments. It belongs clearly to 
$C^{\infty}(P,\tree(I_{k,s}))$.
Finally, since $Q(T)$ is $G_0$-invariant, proposition \ref{commute} implies
that $(L_{h^*}+\gamma(h))Q(T)$ is an
eigenvector of $\cc^{\omega}$ with eigenvalue $\alpha_{k,s}$.
\end{proof}
These rather technical results allow to state the main theorem.
\begin{thm}
If the pair $(V_1,V_2)$ is not critical, then the formula
\[Q_M: (\nabla,T)\mapsto Q_M(\nabla,T)(f)=(p^*)^{-1}[Q_{\omega}(Q(p^*T))(p^*f)],\]
defines a natural and projectively equivariant quantization.
\end{thm}
\begin{proof}
First, the formula makes sense : the function
$Q_{\omega}(Q(p^*T))(p^*f)$ is $H$-equivariant.
The $G_0$-equivariance follows from theorem \ref{hat} and remark~\ref{remf}.
The $\g_1$-equivariance follows from the relations
\[\begin{array}{lll}
L_{h^*}[Q_{\omega}(Q(p^*T))(p^*f)]& =& \L_{h^*}[Q_{\omega}(Q(p^*T))](p^*f)\\
&=&Q_{\omega}[(L_{h^*}+\gamma(h))(Q(p^*T))](p^*f)\\
&=&Q_{\omega}[Q(L_{h^*}(p^*T))](p^*f).
\end{array}\]
Next, the principal symbol of $Q_M(\nabla,T)(f)$ is exactly $T$. Simply notice that
the leading term of $Q(p^*T)$ is $p^*T$, then use the results of
\cite[p. 47]{Capcart}.

Hence, $Q_M(\nabla)$ is a quantization. It is also projectively equivariant by the
very definition of $\omega$.

Eventually, the naturality of the so defined equivariant quantizations is a
consequence of the naturality of all the tools used to derive the formula.
\end{proof}

%%%%
\section{Acknowledgements}
%%%%
F. Radoux thanks the Belgian FRIA for his Research Fellowship.
\bibliographystyle{plain} \bibliography{prquant}

%\noindent Institute of Mathematics,
%\\University of Li\`ege,
%\\Grande Traverse, 12, B 37
%\\B-4000 Sart Tilman
%\\Belgium
%\\mailto:p.mathonet@ulg.ac.be, F.Radoux@ulg.ac.be
\end{document}